\documentclass[12pt]{article}
\usepackage{arxiv}

\usepackage{siunitx}
\usepackage[table]{xcolor}
\definecolor{royal_blue}{rgb}{0.255, 0.412, 1}
\definecolor{medium_blue}{rgb}{0, 0, 0.804}
\usepackage{hyperref}
\hypersetup{colorlinks=true,linkcolor=medium_blue, citecolor=medium_blue}
\usepackage{graphicx}
\usepackage{bm}
\usepackage{amssymb}
\usepackage{multirow}
\usepackage{siunitx}
\usepackage{physics}
\usepackage{subcaption} 
\usepackage{float}
\usepackage{amsmath}
\usepackage{amssymb}
\usepackage{mathdots}
\usepackage{stmaryrd}
\usepackage{threeparttable}
\usepackage{makecell}
\usepackage{tabu,booktabs}

\title{Meshfree implementation of the cardiac monodomain model through the Fragile Points Method}

\author{
  Konstantinos A. ~Mountris\thanks{[mail] kmountris@unizar.es \quad [url] https://www.mountris.org} \\
  Arag\'on Institute of Engineering Research, IIS Arag\'on\\
  CIBER-BBN\\
  University of Zaragoza\\
  Zaragoza, Spain \\
  \texttt{kmountirs@unizar.es} \\
  \And
  Leiting ~Dong \\
  Beihang University, \\
  Beijing, China \\
  \texttt{ltdong@buaa.edu.cn} \\
  \And
  Yue ~Guan \\
  Texas Tech University, \\
  Lubbock, Texas \\
  \texttt{yuguan@ttu.edu} \\
  \And
  Satya N. ~Atluri \\
  Texas Tech University, \\
  Lubbock, Texas \\
  \texttt{snatluri.ttu@gmail.com} \\
  \And
  Esther ~Pueyo \\
  Arag\'on Institute of Engineering Research, IIS Arag\'on, \\
  CIBER-BBN\\
  University of Zaragoza\\
  Zaragoza, Spain \\
  \texttt{epueyo@unizar.es} \\
}

\begin{document}
\maketitle

\begin{abstract}
Meshfree methods for in silico modelling and simulation of cardiac electrophysiology are gaining more and more popularity. These methods do not require a mesh and are more suitable than the Finite Element Method (FEM) to simulate the activity of complex geometrical structures like the human heart. However, challenges such as numerical integration accuracy and time efficiency remain and limit their applicability. Recently, the Fragile Points Method (FPM) has been introduced in the meshfree methods family. It uses local, simple, polynomial, discontinuous functions to construct trial and test functions in the Galerkin weak form. This allows for accurate integration and improved efficiency while enabling the imposition of essential and natural boundary conditions as in the FEM. In this work, we consider the application of FPM for cardiac electrophysiology simulation. We derive the cardiac monodomain model using the FPM formulation and we solve several benchmark problems in 2D and 3D. We show that FPM leads to solutions of similar accuracy and efficiency with FEM while alleviating the need for a mesh. Additionally, FPM demonstrates better convergence than FEM in the considered benchmarks.
\end{abstract}

\keywords{meshfree \and fragile points method \and FPM \and cardiac electrophysiology \and monodomain}

\section{Introduction}
The Finite Element Method (FEM) is widely-used for simulating cardiac electrophysiology mainly due to its robustness and accuracy \cite{heidenreich2010adaptive}. However, accuracy is substantially deteriorated when the mesh undergoes a large deformation or if it does not satisfy specific quality criteria \cite{eppstein2001global}. For this reason, in recent years there has been growing interest in alternative meshfree methods that partially or completely alleviate the mesh-related limitations of FEM. Several meshfree methods have been investigated so far in the context of cardiac electrophysiology. 

The element-free Galerkin (EFG) method offers high convergence rate \cite{belytschko1994element} and it has been used successfully in several applications \cite{belytschko1996meshless,liu2003mesh}. It has been used in cardiac electrophysiology to solve the cardiac monodomain equation using a meshfree representation of the Auckland heart model \cite{zhang2012meshfree}. EFG employs the Moving Least Squares (MLS) approximation \cite{lancaster1981surfaces} for the solution of the Galerkin weak form. Therefore, the imposition of essential boundary conditions requires special treatment due to the lack of the Kronecker delta property in MLS. To solve this issue, Maximum Entropy (MaxEnt) approximants \cite{arroyo2007local,millan2015cell} have been proposed as an alternative to MLS in EFG \cite{mountris2020cell}. Due to the weak Kronecker delta property of MaxEnt, the influence of internal nodes on the boundary of the domain of interest is eliminated and essential boundary conditions can be imposed directly as in FEM. However, the MaxEnt approximation functions are significantly more complex than MLS. As a result, a large number of quadrature points is required to deal with the quadrature accuracy problem of meshfree methods \cite{dolbow1999numerical,babuvska2008quadrature}.

On the other hand, the Mixed Collocation Method (MCM) \cite{mountris2020radial,mountris2021cardiac} is a purely meshfree method that avoids the quadrature accuracy problem. It has been used for the simulation of the cardiac monodomain model demonstrating results in good agreement with FEM \cite{mountris2019novel,mountris2020next}. MCM is a variation of the Meshless Local Petrov-Galerkin method \cite{atluri2006meshless} where the Dirac delta distribution is used as test function and interpolation is applied both on the field function and its gradient \cite{zhang2014meshless,li2008mlpg,atluri2006meshless}. As a result, the integrals in the Petrov-Galerkin weak form are replaced with nodal summation. In addition, essential boundary conditions are imposed directly through collocation. However, being a collocation method, it usually requires the construction of support domains with a large number of collocation points to avoid inaccuracy during the imposition of natural boundary conditions \cite{Libre2008} and may not be as efficient as FEM in large-scale 3D problems.

Very recently, a novel meshfree technique going by the name of Fragile Points Method (FPM) has been added to the artillery of meshfree methods \cite{dong2019new}. FPM uses local, simple, polynomial, discontinuous functions \cite{liszka1980finite} as trial and test functions for the solution of the Galerkin weak form. Using these polynomial trial and test functions, the integration in the Galerkin weak form becomes trivial. Due to the simplicity of the test and trial functions, single-point integration is sufficient to compute integrals with high accuracy in FPM. We refer to table \ref{tab:methods_comp} for a comparison of integration in FPM and other meshfree methods. Moreover, both essential and natural boundary conditions can be imposed as in FEM. However, due to the discontinuity of the trial and test functions, the assembly of the FPM point stiffness matrices leads to an inconsistent global stiffness matrix. To remedy this issue, the numerical flux corrections, which are widely used in Discontinuous Galerkin methods \cite{arnold2002unified}, are employed in FPM to obtain a consistent, sparse and symmetric global stiffness matrix. 

Despite being only recently introduced, FPM has been demonstrated to be an accurate and efficient method with many desired properties (i.e. accurate integration, exact imposition of boundary conditions) and has been already employed to solve linear elasticity \cite{yang2021simple}, heat conduction \cite{guan2020heat,guan2021meshless} and flexoelectric problems with crack propagation \cite{guan2020new}.

In the present study, we employ FPM for the solution of the cardiac monodomain model. Our motivation is to provide a meshfree alternative to FEM for cardiac electrophysiology simulation, while maintaining accuracy and efficiency. The paper is structured as follows. In section \ref{sec:fpm_form}, we describe the theoretical aspects of FPM and the derivation of the cardiac monodomain models. In section \ref{sec:examples}, we present several numerical examples where FPM is applied to simulate electrophysiology in both 2D and 3D benchmark problems, as well as in a large scale biventricular geometry under myocardial infarction conditions. Finally, in section \ref{sec:remarks}, we provide our concluding remarks.

\begin{table}[bt]
	\centering
	\caption{Characteristics of FPM and other numerical methods as in \cite{yang2021simple}.}
	\label{tab:methods_comp}
	\begin{tabu} to \textwidth {cX[l]X[l]X[l]}
		\hline
		\textbf{method} & \textbf{trial function} & \textbf{weak/strong form} & \textbf{numerical integration}\\
		\hline \hline
		FEM & element-based & Galerkin weak form & Gauss integration, inaccurate for highly distorted elements \\
		EFG & point-based continuous & Galerkin weak form & numerical integration with many points \\
		MLPG & point-based continuous & local Petrov-Galerkin weak form & numerical integration with many points in local support domain\\
		FPM & point-based discontinuous & Galerkin weak form with numerical flux corrections & exact integration with one-point-integral for linear trial functions \\
		\hline
	\end{tabu}
\end{table}

\section{Cardiac monodomain model FPM formulation} \label{sec:fpm_form}

\subsection{The cardiac monodomain equation} \label{subsec:monodomain_intro}
We consider the cardiac monodomain model for the simulation of electrical impulse propagation in the human heart, which is governed by the following reaction-diffusion partial differential equation (PDE) \cite{keener2009mathematical}:

\begin{equation} \label{eq:monodomain}
\begin{array}{ll}
    \partial V / \partial t = -I_{ion}(V) / C + \bm{\nabla} \cdot (\bm{D} \bm{\nabla} V) &\textrm{ in } \Omega \\
    \bm{n} \cdot (\bm{D}\bm{\nabla}V) = 0 &\textrm{ in } \partial \Omega
\end{array}
\end{equation}

\noindent where $\partial V / \partial t$ denotes the transmembrane voltage time derivative, $I_{ion}(V)$ the total ionic current, $C$ the cell capacitance per unit surface area and $\bm{D}$ the diffusion tensor. $\Omega$ and $\partial \Omega$ are the domain of interest and its boundary and $\bm{n}$ is the outward unit vector normal to the boundary.

The diffusion tensor $\bm{D}$ is given by:
\begin{equation} \label{eq:diffusion_tensor}
\bm{D} = d_0 \left[(1-\rho)\bm{f}\otimes\bm{f} + \rho \bm{I} \right]
\end{equation}

\noindent where $d_0$ denotes the diffusion coefficient along the cardiac fiber direction, $\bm{f}$ the cardiac fiber direction vector, $\rho \leq 1$ the transverse-to-longitudinal conductivity ratio, $\bm{I}$ the identity matrix and $\otimes$ the tensor product operator.

We employ the operator splitting technique to obtain the decoupled system of Equation (\ref{eq:monodomain}):
\begin{equation} \label{eq:monodomain_decoupled}
\begin{array}{ll}
\partial V(t) / \partial t = -I_{ion}(V(t))/C & \textrm{ in } \Omega \\
\partial V(t) / \partial t = \bm{\nabla} \cdot (\bm{D\nabla}V(t)) & \textrm{ in } \Omega \\
\bm{n} \cdot (\bm{D\nabla}V(t)) = 0 & \textrm{ on } \partial\Omega
\end{array}
\end{equation}

\noindent The decoupled system can be solved efficiently by applying either the Godunov (first-order) or Strang (second-order) methods \cite{schroll2007accuracy}. In the following, we consider the ionic currents $I_{ion}$ of equation (\ref{eq:monodomain_decoupled}) normalized by the capacitance $C$.

\subsection{Computation of trial and test functions in FPM}

FPM is a meshfree method where trial and test functions are established on arbitrarily distributed points in the domain of interest $\Omega$. Unlike mesh-based methods, such as FEM, connectivity information is not required. Simple, local, polynomial trial functions are defined in compact support domains, formed by partitioning the domain in conforming and nonoverlapping subdomains around each point. The partition is not unique and can include subdomains of arbitrary shape, such as polygons that can be obtained by the Voronoi diagram partition.

\begin{figure}[bt]
    \centering
    \includegraphics[width=\textwidth]{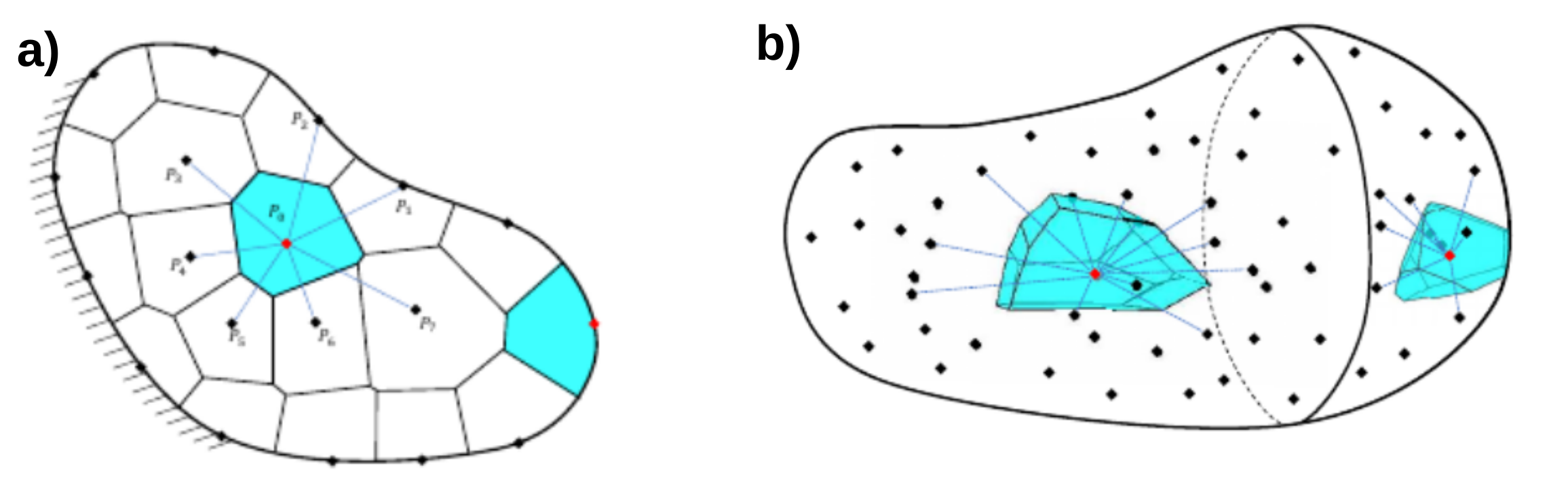}
    \caption{a) Partition of 2D domain with randomly distributed points inside and on its boundary (points $P \in \Omega \cup \partial\Omega$). b) Partition of 3D domain with randomly distributed points inside it (points $P \in \Omega$).}
    \label{fig:fpm_support}
\end{figure}

We define the trial functions at each subdomain in terms of the transmembrane voltage $V$ and its gradient $\nabla V$. The trial function $V_h$ for subdomain $E_0$ that contains the point $P_0$ is given by:
\begin{equation} \label{eq:fpm_trial_func}
    V_h(\bm{x}) = V_0 + (\bm{x}-\bm{x}_0) \cdot \nabla V \Big|_{P_0}, 
\end{equation}

\noindent where $\bm{x}$ is the coordinate vector of a point in $E_0$, $V_0$ is the value of $V_h$ at $P_0$, $\bm{x}_0$ is the coordinate vector of $P_0$ and $\nabla V \Big|_{P_0}$ (i.e. the voltage gradient at $P_0$) is unknown. To compute $\nabla V \Big|_{P_0}$, we employ the Generalized Finite Difference (GFD) method \cite{liszka1980finite} to minimize the weighted discrete L$^2$ norm:
\begin{equation} \label{eq:weighted_l2_norm}
    J = \sum_{i=1}^m w_i \left[(\bm{x}_i - \bm{x}_0) \cdot \nabla V \Big|_{P_0} - (V_i-V_0) \right]^2 \;,
\end{equation}

\noindent where $\bm{x}_i$ denotes the coordinate vector of $P_i \in E_0$, $V_i$ denotes the value of $V_h$ at $P_i$ and $w_i$ denotes the value of the weight function at $P_i$ ($i = 1, 2, \dots ,m$), with $m$ being the number of points in the support domain of $P_0$. We should note that a compact support domain for $P_0$ is defined by the points in the first ring of adjacent subdomains to the subdomain of $P_0$ (Figure \ref{fig:fpm_support}). Assuming constant weight functions, we obtain the transmembrane voltage gradient at $P_0$ by:
\begin{equation} \label{eq:fpm_gradient}
    \nabla V \Big|_{P_0} = \bm{B}\bm{V}_E,
\end{equation}

\noindent where
\begin{equation} \label{eq:gradient_parts}
\begin{split}
\bm{V}_E & = \left[ V_0 \; V_1 \; V_2 \; \dots \; V_m \right]^T, \\
\bm{B} & = (\bm{A}^T\bm{A})^{-1}\bm{A}^T \left[\bm{I}_1 \; \bm{I}_2\right], \\
\bm{I}_1 = \begin{bmatrix}
            -1\\
            -1\\
            \vdots\\
            -1
            \end{bmatrix}_{m\cross1}, \quad \bm{I}_2 = & \begin{bmatrix}
                                                          1 & 0 & \ldots & 0 \\
                                                          0 & 1 & \ddots & \vdots \\
                                                          \vdots & \ddots & \ddots & 0 \\
                                                          0 & \ldots & 0 & 1
                                                         \end{bmatrix}_{m\cross m}, \quad \bm{A} = \begin{bmatrix}
                                                  \bm{x}_1 - \bm{x}_0 \\
                                                  \bm{x}_2 - \bm{x}_0 \\
                                                  \vdots \\
                                                  \bm{x}_m - \bm{x}_0
                                        \end{bmatrix}.
\end{split}
\end{equation}

Using Equation (\ref{eq:fpm_gradient}), the trial function $V_h(\bm{x})$ can be obtained by:
\begin{equation} \label{eq:fpm_shapefun}
    V_h(\bm{x}) = \bm{N}\bm{V}_E, \quad \bm{x} \in E_0
\end{equation}

\noindent where $\bm{N}$ denotes the shape function of $V_h$ in $E_0$:
\begin{equation} \label{eq:fpm_shapefun_explained}
    \bm{N} = [\bm{x}-\bm{x}_0]\bm{B} + \begin{bmatrix}
                                        1 & 0 & \ldots & 0
                                       \end{bmatrix}_{1\times(m+1)}.
\end{equation}

Since the shape function is defined in each subdomain independently, it can be discontinuous at the internal boundaries. The Galerkin weak form in FPM is established by constructing both the trial and test functions using Equation (\ref{eq:fpm_shapefun}). It should be noted that due to the discontinuity of the trial and test functions, the Galerkin weak form will lead to an inconsistent matrix and inaccurate results. Therefore, numerical flux corrections, which are common in Discontinuous Galerkin Finite Element Method \cite{arnold2002unified}, are introduced in FPM to remedy this issue.

\subsection{Numerical flux corrections to remedy the inconsistency}

We address the inconsistency issue by employing the Interior Penalty (IP) numerical flux corrections \cite{mozolevski2007hp}. We start by writing the Galerkin weak form of Equation (\ref{eq:monodomain_decoupled}) for each subdomain $E \in \Omega$, which is given by:
\begin{equation} \label{eq:galerkin_monodomain}
   \int_E v \pdv{V_h}{t} d\Omega + \int_E \nabla v^T \bm{D} \nabla V_h d\Omega = \int_{\partial E} v \bm{n}^T \bm{D} \nabla V_h d\Gamma,
\end{equation}

\noindent where $v$ and $V_h$ are the test and trial functions, respectively, $\partial E$ is the boundary of the subdomain $E$, $\bm{n}$ denotes the outward normal vector to $\partial E$, and $\Gamma$ denotes the set of internal and external boundaries, i.e., $\Gamma = \Gamma_h + \partial \Omega = \Gamma_h + \Gamma_N$, where $\Gamma_h$ is the set of internal boundaries.

Next, the jump operator $\llbracket \; \rrbracket$ and average operator $\{\;\}$ are used to sum Equation (\ref{eq:galerkin_monodomain}) over all subdomains \cite{guan2021meshless}:
\begin{multline} \label{eq:fpm_galerkin_sums}
    \sum_{E\in \Omega} \int_E v \pdv{V_h}{t} d\Omega + \int_E \nabla v^T \bm{D} \nabla V_h d\Omega = \sum_{e\in \Gamma_N} \int_e \llbracket v \rrbracket \{\bm{n}^T\bm{D}\nabla V_h\} d\Gamma +\\
    \sum_{e\in \Gamma_h} \int_e \{v\} \llbracket \bm{n}^T\bm{D}\nabla V_h \rrbracket d\Gamma + \sum_{e\in \Gamma_h} \int_e \llbracket v \rrbracket \{\bm{n}^T\bm{D}\nabla V_h\} d\Gamma,
\end{multline}

\noindent where the jump operator $[\;]$ and average operator $\{\;\}$ are given, $\forall w \in \mathbb{R}$, by:
\begin{equation} \label{eq:jump_avg_operators}
  \llbracket w \rrbracket =
    \begin{cases}
      w |_e^{E1} - w |_e^{E2} & \quad e \in \Gamma_h\\
      w |_e & \quad e \in \partial \Omega
    \end{cases}, \quad
    \{w\} =
    \begin{cases}
      \frac{1}{2} \left(w |_e^{E1} + w |_e^{E2}\right) & \quad e \in \Gamma_h\\
      w |_e & \quad e \in \partial \Omega
    \end{cases}.
\end{equation}

\noindent When $e \in \Gamma_h$ $\left(e \in \partial E_1 \cap \partial E_2\right)$, $\bm{n}$ is a unit vector normal to $e \in \Gamma_h$ and pointing outward from $E$.

Substituting the boundary conditions (Equation \ref{eq:monodomain_decoupled}) into Equation \ref{eq:fpm_galerkin_sums}, we obtain:
\begin{equation}
    \sum_{e\in \Gamma_N} \int_e \llbracket v \rrbracket \{\bm{n}^T\bm{D}\nabla V_h\} d\Gamma = 0.
\end{equation}

Additionally, $\llbracket \bm{n}^T\bm{D}\nabla V_h \rrbracket = 0$ and $\llbracket V_h \rrbracket=0$ when $V_h$ is the exact solution since there is no jump in the internal boundaries. As a result, $\{\bm{n}^T\bm{D}\nabla v\} \llbracket V_h \rrbracket = 0$ and we can replace the term $\{v\} \llbracket \bm{n}^T\bm{D}\nabla V_h \rrbracket$ in Equation (\ref{eq:fpm_galerkin_sums}) with $\{\bm{n}^T\bm{D}\nabla v\} \llbracket V_h \rrbracket$ without affecting the accuracy of the formula.

Finally, the internal penalty numerical flux is applied on $\Gamma_h$ with penalty parameter $\eta$ to obtain the consistent FPM formula \cite{guan2020heat}:
\begin{multline} \label{eq:fpm_consistent}
    \sum_{E\in \Omega} \int_E v \pdv{V_h}{t} d\Omega + \int_E \nabla v^T \bm{D} \nabla V_h d\Omega - \sum_{e\in \Gamma_h} \int_e \{\bm{n}^T\bm{D}\nabla V_h\} \llbracket v \rrbracket d\Gamma - \\ \sum_{e\in \Gamma_h} \int_e \{\bm{n}^T\bm{D}\nabla v\} \llbracket V_h \rrbracket d\Gamma + \sum_{e\in \Gamma_h} \frac{\eta}{h_e} \int_e \llbracket V_h \rrbracket \llbracket v \rrbracket d\Gamma = 0,
\end{multline}

\noindent where $h_e$ is a boundary-dependent parameter with the unit of length. In this work, we define $h_e$ as the distance between points $i$, $j$ when $e\in \partial E_i \cap \partial E_j$. The penalty parameter $\eta$, $\eta>0$, has the same units as $\bm{D}$ and is independent of the boundary size. It should be noted that the penalty parameter should be large enough to ensure stability, but excessively large values should be avoided since they may cause a condition number problem. An extensive discussion on recommended values for the penalty parameter can be found in \cite{guan2020heat}. In this work, we use $\eta$ given by:
\begin{equation} \label{eq:eta}
    \eta = p * \frac{\sum_{i}^{m} \mathcal{V}_{E_i} \Bar{d}_i}{\sum_{i}^{m} \mathcal{V}_{E_i}},
\end{equation}

\noindent where $p$ is a penalty coefficient, $\mathcal{V}_{E_i}$ denotes the volume of the cell containing the $i^{th}$ of the $m$ points in the support domain, and $\Bar{d}_i$ is the mean value of the diagonal entries in the diffusion tensor $\bm{D}_i$ of the $i^{th}$ point.

\subsection{Numerical implementation}
The formula of FPM can be written in matrix form:
\begin{equation} \label{eq:fpm_matrixform}
    \bm{C} \bm{\Dot{V}} + \bm{K}\bm{V} = \bm{0},
\end{equation}

\noindent where $\bm{C}$ and $\bm{K}$ denote the global normalized capacity and diffusion matrices, respectively, and $V$ is the unknown vector collecting the nodal values of the transmembrane potential.

To assemble the global matrices $\bm{C}$ and $\bm{K}$, we substitute the shape function $\bm{N}$ for $V_h$ and $v$ and matrix $\bm{B}$ for $\nabla V_h$ and $\nabla v$ in equation (\ref{eq:fpm_consistent}) to obtain the point normalized capacity matrix $\bm{C}_E$, the point diffusion matrix $\bm{K}_E$, and the internal boundary diffusion matrix $\bm{K}_h$:

\begin{equation} \label{eq:fpm_pointmatrices}
\begin{split} 
    \bm{C}_E =& \int_E \bm{N}^T\bm{N}d\Omega, \quad E\in\Omega, \\
    \bm{K}_E =& \int_E \bm{B}^T\bm{D}\bm{B}d\Omega, \quad E\in\Omega, \\
    \bm{K}_h =& -\frac{1}{2}  \int_e (\bm{N}_1^T\bm{n}_1^T\bm{D}_1\bm{B}_1 + \bm{B}_1^T\bm{D}_1^T\bm{n}_1\bm{N}_1)d\Gamma + \frac{\eta}{h_e}\int_e \bm{N}_1^T\bm{N}_1 d\Gamma \\
    & -\frac{1}{2} \int_e (\bm{N}_2^T\bm{n}_2^T\bm{D}_2\bm{B}_2 + \bm{B}_2^T\bm{D}_2^T\bm{n}_2\bm{N}_2)d\Gamma + \frac{\eta}{h_e}\int_e \bm{N}_2^T\bm{N}_2 d\Gamma \\
    & -\frac{1}{2} \int_e (\bm{N}_1^T\bm{n}_1^T\bm{D}_2\bm{B}_2 + \bm{B}_1^T\bm{D}_1^T\bm{n}_2\bm{N}_2)d\Gamma - \frac{\eta}{h_e}\int_e \bm{N}_1^T\bm{N}_2 d\Gamma \\
    & -\frac{1}{2} \int_e (\bm{N}_2^T\bm{n}_2^T\bm{D}_1\bm{B}_1 + \bm{B}_2^T\bm{D}_2^T\bm{n}_1\bm{N}_1)d\Gamma - \frac{\eta}{h_e}\int_e \bm{N}_2^T\bm{N}_1 d\Gamma, \quad e \in \partial E_1 \cap \partial E_2.
\end{split}
\end{equation}

The assembly of global matrices $\bm{C}$ and $\bm{K}$ is performed as in FEM, where $\bm{C}$ is obtained by assembling all the individual capacity point matrices $\bm{C}_E$, and $\bm{K}$ is obtained by assembling all the point diffusion $\bm{K}_E$ and internal boundary diffusion matrices $\bm{K}_h$. It should be noted that the resulting global diffusion matrix is symmetric, sparse, and positive definite.

\section{Numerical examples} \label{sec:examples}

Simulations for the numerical examples that are presented in the following were performed on a laptop with Intel\textsuperscript{\textregistered} Core\texttrademark i7-4720HQ CPU and 16 GB of RAM. The efficiency of FPM was evaluated by comparing the execution time of FPM simulations with the execution time of FEM using linear elements. The execution time for both FPM and FEM for the considered numerical examples is summarized in Table \ref{tab:times}.

\begin{table}[bt]
\centering
\caption{Required execution time for FPM and FEM solutions for the numerical examples in section \ref{sec:examples}.}
\label{tab:times}
\begin{threeparttable}
\begin{tabular}{cccccc}
\hline
\multicolumn{6}{c}{Electrical propagation in a 2D ventricular tissue (subsection \ref{subsec:ventri2d})} \\
\hline
\hline
\textbf{$\ell$ (mm)} & \textbf{$t_{FEM}$ (min)} & \multicolumn{4}{c}{\textbf{$t_{FPM}$ (min)}}\\
     & & $p=1$ & $p=2$ & $p=5$ & $p=10$ \\
\hline
2 & 0.87 & 0.90 & 0.95 & 0.88 & 0.84\\
1 & 1.52 & 1.59 & 1.63 & 1.51 & 1.51  \\
0.5 & 3.7 & 3.65 & 3.77 & 3.76 & 4.21 \\
0.25 & 12.33 & 12.46 & 12.47 & 14.98 & 17.26 \\
\hline
\multicolumn{6}{c}{Acetylcholine-induced effects in human atrial electrical activity (subsection \ref{subsec:atria2d})} \\
\hline \hline
\multicolumn{2}{c}{\textbf{$ACh$ ($\mu$M)}} & \multicolumn{2}{c}{\textbf{$t_{FEM}$ (min)}} & \multicolumn{2}{c}{\textbf{$t_{FPM}$ (min)}}\\
\hline
\multicolumn{2}{c}{0.00} & \multicolumn{2}{c}{13.04} & \multicolumn{2}{c}{15.21} \\
\multicolumn{2}{c}{0.02} & \multicolumn{2}{c}{13.82} & \multicolumn{2}{c}{15.49} \\
\multicolumn{2}{c}{0.04} & \multicolumn{2}{c}{13.26} & \multicolumn{2}{c}{15.42} \\
\multicolumn{2}{c}{0.06} & \multicolumn{2}{c}{13.62} & \multicolumn{2}{c}{14.67} \\
\multicolumn{2}{c}{0.08} & \multicolumn{2}{c}{14.62} & \multicolumn{2}{c}{15.93} \\
\multicolumn{2}{c}{0.10} & \multicolumn{2}{c}{13.26} & \multicolumn{2}{c}{14.61} \\
\hline
\multicolumn{6}{c}{Electrical propagation in a benchmark 3D cuboid geometry (subsection \ref{subsec:cuboid})} \\
\hline \hline
\multicolumn{2}{c}{\textbf{$\ell$ (mm)}} & \multicolumn{2}{c}{\textbf{$t_{FEM}$ (min)}} & \multicolumn{2}{c}{\textbf{$t_{FPM}$ (min)}}\\
\hline
\multicolumn{2}{c}{0.5} & \multicolumn{2}{c}{0.35} & \multicolumn{2}{c}{0.38} \\
\multicolumn{2}{c}{0.2} & \multicolumn{2}{c}{2.94} & \multicolumn{2}{c}{3.26} \\
\multicolumn{2}{c}{0.1} & \multicolumn{2}{c}{27.9} & \multicolumn{2}{c}{29.25} \\
\hline
\multicolumn{6}{c}{Simulation of electrical activation in 3D biventricular infarction model (subsection \ref{subsec:pacing})} \\
\hline \hline
\multicolumn{2}{c}{\textbf{pacing}} & \multicolumn{2}{c}{\textbf{$t_{FEM}$ (min)}} & \multicolumn{2}{c}{\textbf{$t_{FPM}$ (min)}}\\
\multicolumn{2}{c}{basal} & \multicolumn{2}{c}{9.37} & \multicolumn{2}{c}{11.52} \\
\multicolumn{2}{c}{apical} & \multicolumn{2}{c}{9.15} & \multicolumn{2}{c}{11.46} \\
\hline
\end{tabular}
\end{threeparttable}
\end{table}

\subsection{Electrical propagation in a 2D ventricular tissue} \label{subsec:ventri2d}

We considered a 4x4 cm human ventricular tissue with fibers aligned parallel to the x-axis. The longitudinal diffusion coefficient was $d_0 = 0.0013$ cm$^2$/ms and the transverse-to-longitudinal conductivity ratio was $\rho = 0.15$. Electrophysiology was modeled by the O'Hara et al. \cite{ohara2011simulation} human ventricle action potential (AP) model for epicardial cells. Periodic stimuli of duration $t_d = 1$ ms and amplitude $A$ equal to twice the diastolic threshold were applied on the left side of the tissue ($x = 0$ cm) at a frequency $f = 1$ Hz. AP propagation was simulated for a total time $t = 3$ s, after achievement of steady-state, using the dual adaptive explicit time integration method (DAETI) \cite{mountris2021dual} with time step $dt = 0.1$ ms.

\begin{figure}[bt]
    \centering
    \includegraphics[width=\textwidth]{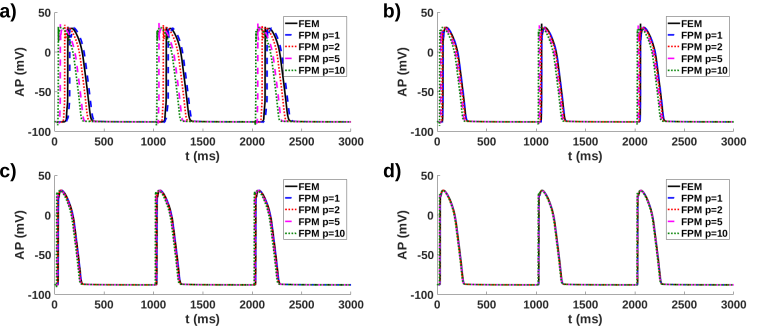}
    \caption{AP propagation in 2D human ventricular tissue for FEM (black), FPM with penalty coefficient $p=1$ (blue dashed), $p=2$ (red dotted), $p=5$ (magenta dashed) and $p=10$ (green dotted). The four panels correspond to nodal discretizations with spacing $\ell$ of 0.25 mm (a), 0.5 mm (b), 1 mm (c) and 2 mm (d).}
    \label{fig:ventri_ap_comp}
\end{figure}

Solutions obtained by FPM with different penalty coefficients, $p = \{1,2,5,10\}$, were compared to a solution obtained by FEM with bilinear isoparametric elements (Figure \ref{fig:ventri_ap_comp}). We considered four different nodal discretizations with spacing $\ell = \{0.25, 0.5, 1, 2\}$ mm, each of which is presented in one of the panels of Figure \ref{fig:ventri_ap_comp}. Differences between FEM and FPM solutions were evaluated in terms of conduction velocity ($CV$) and AP duration (APD) at 90\% repolarization. 

The maximum percentage difference between FPM and FEM in terms of $CV$ was 116.9\%, while in terms of $APD_{90}$ it was 2.7\% for nodal spacing $\ell = 2$ mm and $p = 10$. The minimum percentage difference between FPM and FEM in terms of $CV$ was 7.1\%, while in terms of $APD_{90}$ it was 0.0\% for nodal spacing $\ell = 0.25$ mm and $p = 1$. For a given nodal discretization, the best agreement between FEM and FPM solutions was always obtained for $p = 1$. For penalty factor values $p = 5$ and $p = 10$, the $CV$ values obtained by FPM with coarse nodal discretizations $\ell = 2$ mm and $\ell = 1$ mm were in better agreement with the $CV$ value obtained by FEM for a dense nodal discretization ($\ell = 0.25$ mm). 

\begin{figure}[bt]
    \centering
    \includegraphics[width=\textwidth]{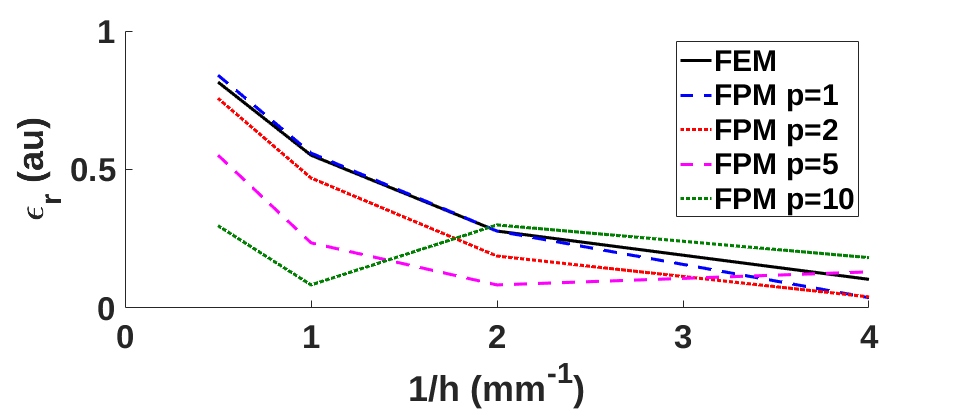}
    \caption{Convergence analysis in terms of $CV$ for FEM (black), FPM with penalty coefficient $p=1$ (blue dashed), $p=2$ (red dotted), $p=5$ (magenta dashed) and $p=10$ (green dotted).}
    \label{fig:conv}
\end{figure}

To evaluate the effect of the penalty coefficient, we performed a convergence analysis where we evaluated the relative error ($\epsilon_r$) in $CV$ for the solutions obtained by FPM with $p = \{1,2,5,10\}$ with respect to the FEM solution for nodal spacing $\ell = \{0.25, 0.5, 1, 2\}$ mm:
\begin{equation} \label{eq:rel_error}
    \epsilon_r = \frac{CV_a-CV_e}{CV_e},
\end{equation}

\noindent where $CV_a$ is the $CV$ value obtained by FPM or FEM for each of the tested nodal discretizations and $CV_e$ is the value obtained by FEM for a dense nodal discretization with $\ell = 0.1$ mm. The relative error, $\epsilon_r$, for FPM with $p=10$ and $p=5$ was smaller than for FEM for coarse discretizations, but it was larger for finer discretizations (Figure \ref{fig:conv}). For FPM with $p=2$ and $p=1$, $\epsilon_r$ decreased monotonically and its values were smaller than for FEM for practically all nodal spacing values. In the following numerical examples, we used $p = [1,2]$.  

\subsection{Acetylcholine-induced effects in human atrial electrical activity} \label{subsec:atria2d}

Acetylcholine-induced APD shortening in atrial myocytes facilitates the initation and perpetuation of atrial fibrillation \cite{bayer2019acetylcholine,smeets1986wavelength}, which is a major risk factor for ischemic stroke. The parasympathetic neurotransmitter acetylcholine (ACh) shortens APD by activating the ACh-sensitive inward rectifier potassium current, I$_{KACh}$. Simulation of ACh-induced alterations in atrial electrophysiology is of interest in atrial fibrillation research \cite{celotto2021location}.

Here, we considered a 2D atrial tissue of $4\times4$ cm with fibers aligned parallel to the X-axis in which ACh was distributed homogeneously throughout the tissue. The atrial myocyte electrophysiology was described using the Maleckar et al. model \cite{maleckar2008mathematical}. The longitudinal diffusion coefficient was set to $d_0 = 0.0035$ cm$^2$/ms and the transverse-to-longitudinal conductivity ratio was set to $\rho = 0.5$. The same periodic stimulation protocol as in subsection \ref{subsec:ventri2d} was simulated for $t_s = 3$ s and the DAETI method with time step $dt=0.05$ ms was used to numerically solve AP propagation. 

\begin{figure}[bt]
    \centering
    \includegraphics[width=\textwidth]{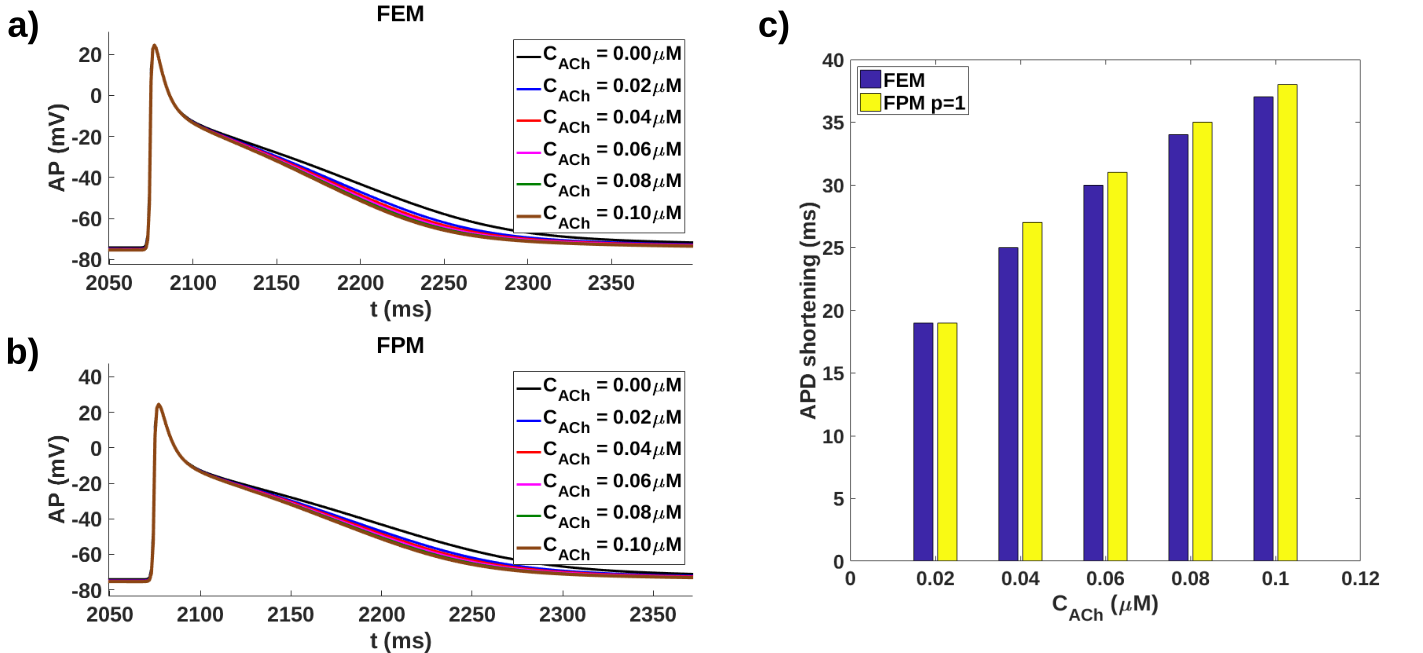}
    \caption{APs at the center of a $4\times4$ cm atrial tissue for different ACh concentrations obtained from a) FEM and b) FPM simulations. c) APD shortening induced by acetylcholine concentration for FEM (blue) and FPM (yellow) simulations.}
    \label{fig:apd_short}
\end{figure}

Simulations were performed using FPM with penalty coefficient $p = 1$ considering a regular nodal distribution with $\ell = 0.25$ mm and were compared with FEM simulations using bilinear isoparametric elements. The induced APD shortening was recorded at the center of the tissue for different ACh concentrations, $C_{ACh} = \{0.02, 0.04, 0.06, 0.08, 0.1\}$ $\mu$M (Figure \ref{fig:apd_short}). APD shortening was found to range from 19 ms for the lowest ACh concentration to 38 ms for the highest ACh concentration. FPM results were found in good agreement with FEM, with differences of up to 2 ms. These results are in line with previous studies reporting APD shortening of up to 40 ms for $C_ACh = 0.1$ $\mu$M \cite{bayer2019acetylcholine,celotto2019calcium,celotto2020sk}.  

\subsection{Electrical propagation in a benchmark 3D cuboid geometry} \label{subsec:cuboid}

\begin{figure}[bt]
\centering
\includegraphics[width=\textwidth]{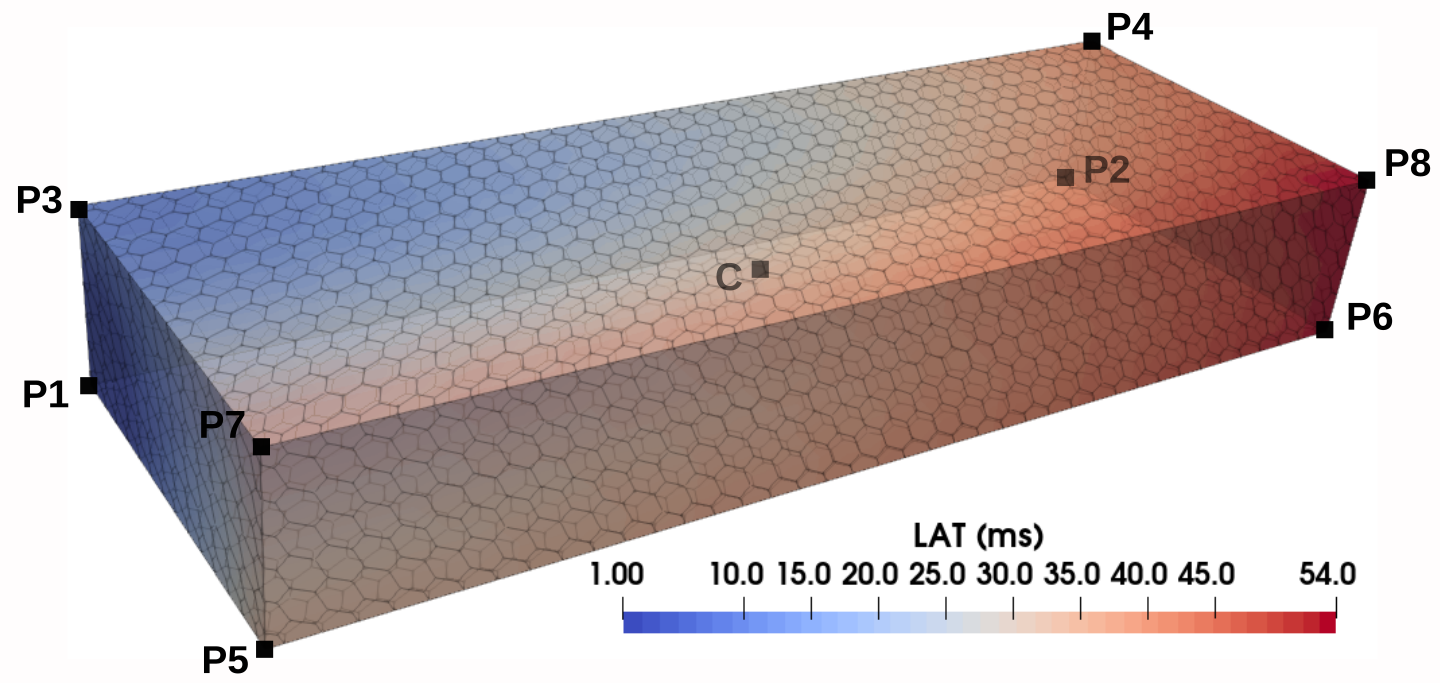}
\caption{Local activation time (LAT) map for the 3D cuboid benchmark geometry described in \cite{niederer2011verification} with space discretization $\ell=0.5$ mm.}
\label{fig:bench}
\end{figure}

In this example, we solved the standard benchmark problem for verification of cardiac tissue electrophysiology simulators described in \cite{niederer2011verification}. The problem considered the electrical stimulation of a 3D cuboid of human ventricular tissue with dimensions $3 \times 7 \times 20$ mm and cardiac fibers parallel to the Z axis. Epicardial cell electrophysiology was described by the Ten Tusscher et al. model \cite{ten2006alternans}. The longitudinal diffusion coefficient was set to $d_0=0.00115$ cm$^2$/ms and the transverse-to-longitudinal ratio was set to $\rho = 0.12$. A periodic stimulus with frequency $f = 1$ Hz, amplitude $A = 50 mA$ and duration $t_d = 2$ ms was delivered at a cubic region with dimensions $1.5 \times 1.5 \times 1.5$ mm located at corner P1 (Figure \ref{fig:bench}).

The activation time at the corners (P1 -- P8) and the center (C) of the cuboid were recorded for three different nodal discretizations, $\ell = \{0.1, 0.2, 0.5\}$ mm, for a simulation using FPM with penalty coefficient $p = 1$. Integration was performed using the DAETI method with time step $dt = 0.1$ ms. The activation times obtained by FPM were compared with those of FEM (Table \ref{tab:bench}) previously reported in \cite{mountris2021dual} and validated with the reported activation times in \cite{niederer2011verification}.

The activation times for FPM and FEM for nodal spacing $\ell = 0.1$ mm were found to be in good agreement. For coarse nodal discretizations, the FEM solution led to larger activation times as compared to FPM, especially at points \textbf{P5}--\textbf{P8}. Activation times at \textbf{P7} for FEM simulations were found to be $\times2.37$ larger for the spacing $\ell = 0.5$ mm than for the spacing $\ell = 0.1$ mm. On the other hand, this difference was found to be $\times1.5$ for FPM simulations. For a nodal spacing $\ell = 0.2$ mm, activation times were notably closer to those of $\ell = 0.1$ mm when using FPM than when using FEM. These results demonstrate the higher convergence of FPM compared to FEM and the improved accuracy for coarse nodal discretizations, in good agreement with the findings in subsection \ref{subsec:ventri2d}.

\begin{table}[bt]
\centering
\caption{Activation times at the corners (P1--P8) and the center (C) of the 3D cuboid benchmark problem.}
\label{tab:bench}
\begin{threeparttable}
\begin{tabular}{cccccccccc}
\hline
\textbf{h (mm)} & \textbf{P1} & \textbf{P2} & \textbf{P3} & \textbf{P4} & \textbf{P5} & \textbf{P6} & \textbf{P7} & \textbf{P8} & \textbf{C} \\
\hline
\multicolumn{10}{c}{FPM activation time (ms)} \\
\hline
0.5 & 1 & 40 & 10 & 40 & 33 & 53 & 36 & 54 & 25 \\
0.2 & 1 & 30 & 8 & 31 & 22 & 39 & 25 & 40 & 19 \\
0.1 & 1 & 29 & 7 & 30 & 22 & 37 & 24 & 38 & 18 \\
\hline
\multicolumn{10}{c}{FEM activation time (ms) as in \cite{mountris2021dual}} \\
\hline
0.5 & 1 & 49 & 22 & 58 & 94 & 109 & 98 & 111 & 54 \\
0.2 & 1 & 31 & 11 & 35 & 35 & 51 & 39 & 54 & 25 \\
0.1 & 1 & 29 & 8 & 31 & 27 & 41 & 29 & 43 & 20 \\
\hline
\end{tabular}
\end{threeparttable}
\end{table}

\subsection{Simulation of electrical activation in a 3D biventricular infarction model} \label{subsec:pacing}

We employed the FPM method to simulate AP propagation in a 3D biventricular model under myocardial infarction conditions and we compared the activation pattern obtained by FPM with the one obtained by FEM. We considered a simplified model where zero conduction was assumed for the scar tissue, while border zone effects were not taken into account. 

\begin{figure}[bt]
    \centering
    \includegraphics[width=0.6\textwidth]{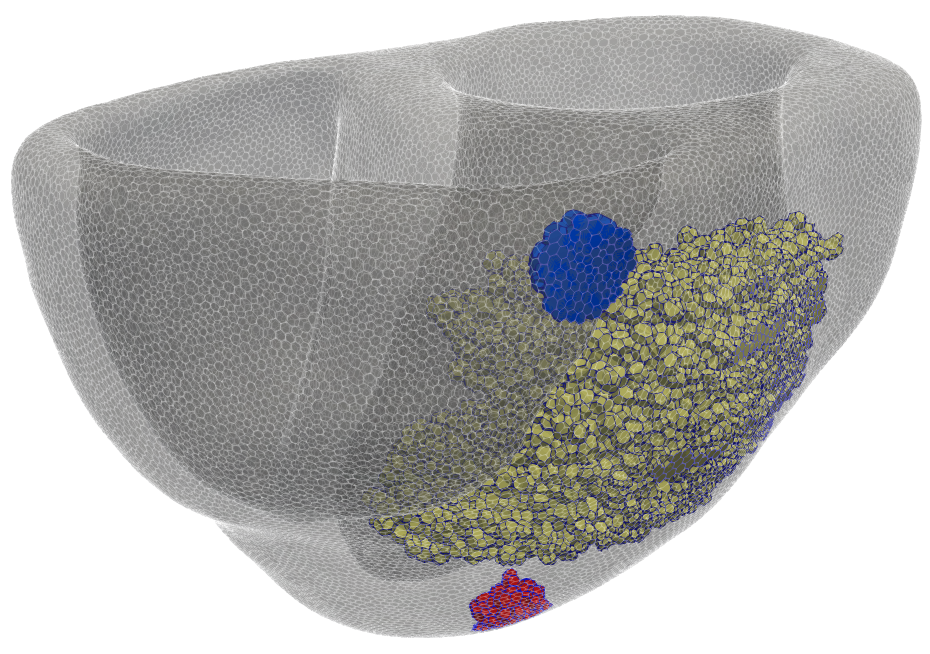}
    \caption{Biventricular model depicting the infarcted region (yellow), the basal pacing region (blue) and the apical pacing region (red).}
    \label{fig:biventri_geo}
\end{figure}

The biventricular anatomy was constructed using ex vivo diffusion weighted imaging (DWI) of a porcine heart. The two ventricles of the heart were manually segmented and a scar was introduced at the septal-anterior wall of the left ventricle (Figure \ref{fig:biventri_geo}). A tetrahedral mesh (nodes: 70521, elements: 311150) was generated from the segmented data using iso2mesh \cite{fang2009tetrahedral}. FPM cells were generated for the nodes of the tetrahedral mesh using a dual polyhedral mesh generation algorithm \cite{garimella2014polyhedral,kim2014efficient}. Fiber direction was determined by computing the diffusion tensors using an algorithm based on Riemannian distances \cite{barmpoutis2010unified}.

Cell electrophysiology was modeled by the O'Hara et al. model \cite{ohara2011simulation} considering endocardium:midmyocardium:epicardium with 50:20:30 ratio. Pacing was performed by applying a periodic stimulus with amplitude equal to twice the diastolic threshold and frequency $f = 1$ Hz at two tested stimulation regions: one at the base of the model located at the anterior wall of the right ventricular base near the septum and the other one at the apex of the left ventricle (Figure \ref{fig:biventri_geo}). 

\begin{figure}[bt]
    \centering
    \includegraphics[width=\textwidth]{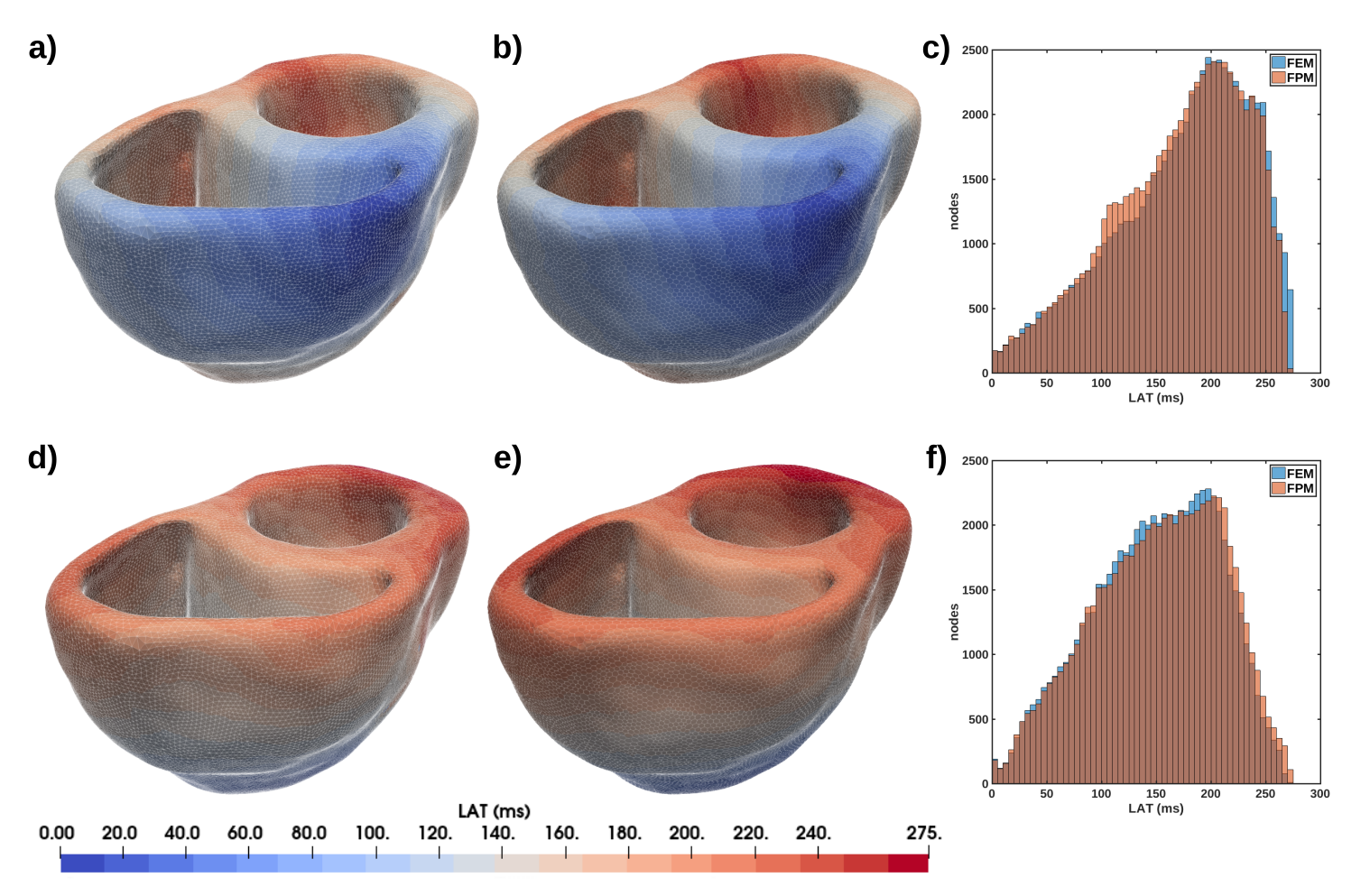}
    \caption{LAT for a) basal pacing using FEM, b) basal pacing using FPM with $p=1.5$, d) apical pacing using FEM, e) apical pacing using FPM with $p=1.5$.  Comparisons of LAT histograms are shown for basal pacing in c)  and for apical pacing in f).}
    \label{fig:biventri_lat}
\end{figure}

The local activation time (LAT) at each node of the model was computed using FPM with penalty coefficient $p = 1.5$ and it was compared with the LAT value obtained from a FEM simulation. Mean LAT for basal pacing was 170 ms for FPM and 173 ms for FEM, while mean LAT for apical pacing was 151 ms for FPM and 148 ms for FEM. The LAT histograms presented in figure \ref{fig:biventri_lat} demonstrate the good agreement between FPM and FEM simulations, with FPM rendering a valuable alternative to FEM for large scale cardiac electrophysiology simulations.

\section{Concluding remarks} \label{sec:remarks}

In this work, we presented the Fragile Points Method (FPM) for in silico cardiac electrophysiology applications. Despite being common for meshfree methods to sacrifice efficiency for accuracy and vice versa, FPM was proven to achieve similar accuracy and efficiency to FEM (Table \ref{tab:times}). 

We found that solutions obtained by FPM were in good agreement with FEM for both 2D and 3D scenarios, allowing to simulate action potential propagation with high accuracy under different physiological and pathological conditions. FPM demonstrated higher convergence than FEM (Table \ref{tab:bench}), in line with previous findings. This may have an important role in large scale applications where FPM could increase time efficiency without losing accuracy. By adjusting the penalty coefficient improved solutions could be obtained even for coarse discretizations. However, for large values of the penalty coefficient ($p$) accuracy could be deteriorated for fine discretizations. We found $p = [1,2]$ to lead to accurate results for all discretizations.

Finally, the ability of FPM to provide similar accuracy and efficiency to FEM without requiring mesh connectivity information renders the method an interesting alternative to FEM, particularly for personalized image-based modeling applications.

\section*{Acknowledgements}

This work was supported by MCIN/AEI/10.13039/501100011033 (Spain) through project PID2019-105674RB-I00, by the European Research Council under grant agreement ERC-StG 638284, by the European Union’s H2020 Program under grant agreement No. 874827 (BRAV3) and by European Social Fund (EU) and Arag\'on Government through BSICoS group (T39\_20R). Computations were performed by the ICTS NANBIOSIS (HPC Unit at University of Zaragoza).

\bibliographystyle{unsrt}  
\bibliography{main}

\end{document}